\documentclass[
a4paper,
11pt,
twoside,
]{article}
\usepackage[utf8]{inputenc}
\usepackage[T1]{fontenc}
\usepackage[osf,sc,noBBpl]{mathpazo} 
\usepackage[scaled]{beramono}
\usepackage{microtype}
\usepackage[english]{babel}
\usepackage{geometry}
\usepackage{amsmath,amsfonts,amssymb,amsthm}
\usepackage{etoolbox}
\usepackage{graphicx}
\usepackage[dvipsnames]{xcolor}
\usepackage{mathtools}
\usepackage{siunitx}
\usepackage{booktabs}
\usepackage{enumitem}
\setlist{itemsep=0pt}

\newcommand{\D}{\mathop{}\!\mathrm{d}}
\DeclareMathOperator{\diag}{diag}

\newcommand{\E}{\mathrm{e}}
\renewcommand{\restriction}{\raise-.5ex\hbox{\ensuremath{\upharpoonright}}}

\usepackage{fnpct}

\definecolor{webgreen}{rgb}{0,.5,0}
\definecolor{webbrown}{rgb}{.6,0,0}
\definecolor{myblue}{rgb}{0,0.25,0.5}
\usepackage[
  colorlinks=true,
  breaklinks=true,
  bookmarksnumbered,
  pdfhighlight=/O,
  urlcolor=webbrown,
  linkcolor=myblue,
  citecolor=webgreen,
  hyperfootnotes=true,
]{hyperref}
\newcommand{\email}[1]{\href{mailto:#1}{\texttt{#1}}}
\makeatletter
\newcommand{\pagerefstar}{\@pagerefstar}
\makeatother

\usepackage{csquotes}

\usepackage[
  mincrossrefs=99,
  style=numeric-comp,
  bibstyle=siam,
  sorting=nyt,
  backref=true,
  backend=biber,
  url=true,
  doi=true,
  maxbibnames=99,
  maxnames=99,
  giveninits=true,
  useprefix=true,
  date=year,
  autolang=other
]{biblatex}
\DeclareCiteCommand{\citeauthornsc}
  {%
   \boolfalse{citetracker}%
   \boolfalse{pagetracker}%
   \usebibmacro{prenote}}
  {\ifciteindex
     {\indexnames{labelname}}
     {}%
   \printnames{labelname}}
  {\multicitedelim}
  {\usebibmacro{postnote}}
\DeclareSourcemap{
  \maps[datatype=bibtex]{
    \map{
      \step[fieldsource=doi,final]
      \step[fieldset=url,null]
    }
  }
}

\usepackage[capitalize,nameinlink,sort]{cleveref}[0.19]
\crefname{section}{section}{sections}
\crefname{subsection}{subsection}{subsections}
\Crefname{figure}{Figure}{Figures}
\crefformat{equation}{\textup{#2(#1)#3}}
\crefrangeformat{equation}{\textup{#3(#1)#4--#5(#2)#6}}
\crefmultiformat{equation}{\textup{#2(#1)#3}}{ and \textup{#2(#1)#3}}
{, \textup{#2(#1)#3}}{, and \textup{#2(#1)#3}}
\crefrangemultiformat{equation}{\textup{#3(#1)#4--#5(#2)#6}}%
{ and \textup{#3(#1)#4--#5(#2)#6}}{, \textup{#3(#1)#4--#5(#2)#6}}{, and \textup{#3(#1)#4--#5(#2)#6}}
\Crefformat{equation}{#2Equation~\textup{(#1)}#3}
\Crefrangeformat{equation}{Equations~\textup{#3(#1)#4--#5(#2)#6}}
\Crefmultiformat{equation}{Equations~\textup{#2(#1)#3}}{ and \textup{#2(#1)#3}}
{, \textup{#2(#1)#3}}{, and \textup{#2(#1)#3}}
\Crefrangemultiformat{equation}{Equations~\textup{#3(#1)#4--#5(#2)#6}}%
{ and \textup{#3(#1)#4--#5(#2)#6}}{, \textup{#3(#1)#4--#5(#2)#6}}{, and \textup{#3(#1)#4--#5(#2)#6}}
\crefdefaultlabelformat{#2\textup{#1}#3}

\crefname{chapter}{Chapter}{Chapters}
\crefname{appendix}{Appendix}{Appendices}
\crefname{subappendix}{Section}{Sections}
\Crefname{subappendix}{Section}{Sections}
\crefname{page}{page}{pages}
\crefname{section}{Section}{Sections}
\Crefname{section}{Section}{Sections}
\crefname{subsection}{Section}{Sections}
\Crefname{subsection}{Section}{Sections}
\crefname{remark}{Remark}{Remarks}
\crefname{footnote}{Footnote}{Footnotes}

\theoremstyle{plain}
\newtheorem{remark}{Remark}[section]

\usepackage{caption}
\usepackage{fancyhdr}
\usepackage{lastpage}

\title{A practical approach to computing Lyapunov exponents of renewal and delay equations}

\author{Dimitri Breda, Davide Liessi}

\newcommand{\mydate}{29 December 2023}

\fancypagestyle{mypagestyle}{%
  \fancyhf{}%
  \fancyhead[LO,RE]{\scriptsize A practical approach to computing Lyapunov exponents of renewal and delay equations}%
  \fancyhead[RO,LE]{\scriptsize Breda, Liessi}%
  \fancyfoot[LE,RO]{\scriptsize \thepage\,/\,\pagerefstar{LastPage}}%
  \fancyfoot[LO,RE]{\scriptsize \mydate{}}%
}
\begin{document}

\thispagestyle{empty}

\begin{center}
\LARGE
A practical approach to computing Lyapunov exponents of renewal and delay equations

\bigskip
\large
Dimitri Breda\footnote{\email{dimitri.breda@uniud.it}},
Davide Liessi\footnote{\email{davide.liessi@uniud.it}}

\medskip
\small
CDLab -- Computational Dynamics Laboratory \\
Department of Mathematics, Computer Science and Physics, University of Udine \\
Via delle Scienze 206, 33100 Udine, Italy

\medskip
\large
\mydate
\end{center}

\begin{abstract}
We propose a method for computing the Lyapunov exponents of renewal equations (delay equations of Volterra type) and of coupled systems of renewal and delay differential equations.
The method consists of the reformulation of the delay equation as an abstract differential equation, the reduction of the latter to a system of ordinary differential equations via pseudospectral collocation and the application of the standard discrete QR method.
The effectiveness of the method is shown experimentally and a MATLAB implementation is provided.

\end{abstract}

\renewcommand{\thefootnote}{\arabic{footnote}}

\section{Introduction}
\label{sec:intro}

A delay equation is a functional equation consisting of
``a rule for extending a function of time towards the future on the basis of the (assumed to be) known past'' \cite{DiekmannVerduynLunel2021}.
A renewal equation (RE) is a delay equation of Volterra type, i.e., the rule for extension prescribes the value of the unknown function itself, instead of the value of its derivative, as in the case of delay differential equations (DDEs).

The goal of this work is to compute the (dominant) Lyapunov exponents (LEs) of REs and of coupled systems of REs and DDEs (henceforth \emph{coupled} equations).
The usefulness of LEs for measuring the asymptotic exponential behavior of solutions is well known; for example, they can be used to study the average asymptotic stability of solutions, the insurgence of chaotic dynamics and the effects of perturbations on the system, as well as to estimate the entropy or the dimension of attractors.

As for DDEs, recent methods for computing the LEs have been proposed, particularly \cite{BredaVanVleck2014} and \cite{BredaDellaSchiava2018}, which use two different approaches (for other methods, see the references in the cited works).

In \cite{BredaDellaSchiava2018} the DDE is reformulated as an abstract differential equation and a pseudospectral discretization is applied \cite{BredaDiekmannGyllenbergScarabelVermiglio2016}, yielding a system of ordinary differential equations (ODEs); LEs are then computed by using the standard discrete QR method (henceforth \emph{DQR}) for ODEs proposed in \cite{DieciVanVleck2004,DieciJollyVanVleck2011}.
In \cite{BredaVanVleck2014}, instead, the problem is tackled directly: the DDE is posed in an infinite-dimensional Hilbert space as the state space, the associated family of evolution operators is discretized and the DQR is adapted and applied to the finite-dimensional approximation; for the error analysis, the DQR is raised to infinite dimension and compared to the approximated DQR used for the computations.

As for REs, as far as we know, there are no methods available in the literature for computing LEs.
Only a first example of naive computation can be found in \cite{BredaDiekmannLiessiScarabel2016}, where it is done simply to exemplify the versatility of the collocation techniques used therein, without attempting a rigorous formulation and error analysis.

In the present work of computational nature, we develop a practical method following the approach of \cite{BredaDellaSchiava2018} described above, and based on \cite{ScarabelDiekmannVermiglio2021} for the reformulation of REs into abstract differential equations.
As in \cite{BredaDellaSchiava2018}, we use the DQR to compute the LEs of the approximating ODE, but, in principle, any method can be used; our choice was motivated by our goal of providing a practical way of computing LEs by using ready-to-use code for a well-known method.

In \cref{sec:dqr-ode} we recall the DQR for linear ODEs.
Then, in \cref{sec:pscoll} we define the reformulation of REs, DDEs and coupled equations into abstract differential equations and their pseudospectral collocation into ODEs.
After describing the implementation choices in \cref{sec:impl}, we present in \cref{sec:res} some numerical experiments concerning the convergence of the method for an example RE with many known properties, as well as some examples of computation of LEs of REs and coupled equations.
Finally, we present some concluding remarks in \cref{sec:concluding}.

The MATLAB codes implementing the method and the scripts to reproduce the experiments of \cref{sec:res} are available at \url{http://cdlab.uniud.it/software}.

\section{DQR for ODEs}
\label{sec:dqr-ode}

In this \lcnamecref{sec:dqr-ode}, we first illustrate the DQR to compute the LEs of linear nonautonomous ODEs; in the nonlinear case, one previously linearizes around a reference trajectory in the attractor.
Then, we comment on the relevant literature.

\bigskip

Let $n$ be a positive integer and consider the ODE
\begin{equation}\label{ODE}
z'(t)=A(t)z(t)
\end{equation}
for $A\colon[0,+\infty)\to\mathbb{R}^{n\times n}$ continuous and bounded; also, let $Z(t)$ be the fundamental matrix solution exiting from a given nonsingular matrix $Z_{0}\in\mathbb{R}^{n\times n}$ prescribed at time $0$ without loss of generality.
For any sequence $\{t_{k}\}_{k\in\mathbb{N}}$ of time instants strictly increasing from $t_{0}=0$, construct the iterative QR factorization%
\footnote{In what follows, a QR factorization of a nonsingular matrix is intended as the unique one with positive diagonal elements.}
\begin{equation}\label{ZQRk1}
Z(t_{k})=Q_{k}R_{k}
\end{equation}
starting from $Z_{0}=Q_{0}R_{0}$ and, at each step $j=1,\dots,k$, solving the $n$ initial value problems (IVPs)
\begin{equation}\label{IVPGAMMA}
\left\{
\begin{aligned}
\Gamma'(t,t_{j-1}) &= A(t)\Gamma(t,t_{j-1}),\quad t\in[t_{j-1},t_{j}],\\
\Gamma(t_{j-1},t_{j-1}) &= Q_{j-1}
\end{aligned}
\right.
\end{equation}
and factorizing the solution at $t_{j}$ as
\begin{equation}\label{GAMMAQR}
\Gamma(t_{j},t_{j-1})=Q_{j}R_{j,j-1}.
\end{equation}
If $S(t,s)\coloneqq Z(t)Z(s)^{-1}$ is the state transition matrix associated with \cref{ODE}, then
\begin{equation*}
\begin{split}
Z(t_{k})={}&S(t_{k},t_{k-1})\cdots S(t_{2},t_{1})S(t_{1},t_{0})Q_{0}R_{0}\\
={}&S(t_{k},t_{k-1})\cdots S(t_{2},t_{1})\Gamma(t_{1},t_{0})R_{0}\\
={}&S(t_{k},t_{k-1})\cdots S(t_{2},t_{1})Q_{1}R_{1,0}R_{0}\\
={}&S(t_{k},t_{k-1})\cdots\Gamma(t_{2},t_{1})R_{1,0}R_{0}\\
={}&S(t_{k},t_{k-1})\cdots Q_{2}R_{2,1}R_{1,0}R_{0}\\
\cdots\\
={}&Q_{k}R_{k,k-1}\cdots R_{1,0}R_{0}.
\end{split}
\end{equation*}
The uniqueness of the QR factorization and \cref{ZQRk1} give
\begin{equation*}
R_{k}=\left(\prod_{j=1}^{k}R_{j,j-1}\right)R_{0},
\end{equation*}
so that, eventually, (upper%
\footnote{%
Lower exponents come either as $\liminf$ or as upper exponents of the adjoint system.
Note, however, that for \emph{regular} ODEs in the sense of Lyapunov (see, e.g., \cite[Definition 3.5.1]{Adrianova1995}) the LEs exist as exact limits, and such quantities are meaningful for stability statements in the original nonlinear system, thus avoiding the so-called Perron effect (see, e.g., \cite{LeonovKuznetsov2007}).
Nevertheless, an in-depth study of the theoretical requirements of the original delay equation guaranteeing the regularity of the ODE obtained by pseudospectral collocation is beyond the scope of the present work, as well as the extension of the regularity concept itself to the infinite-dimensional case of delay equations (extension that the authors, to the best of their knowledge, are not aware of; see, anyway, \cite{BarreiraValls2017} and the references therein).
}%
) LEs are recovered as
\begin{equation}\label{LE}
\lambda_{i}=\limsup_{k\to\infty}\frac{1}{t_{k}}\sum_{j=1}^{k}\ln[R_{j,j-1}]_{i,i},\quad i=1,\dots,n.
\end{equation}
Above, $[R_{j,j-1}]_{i,i}$ denotes the $i$-th diagonal entry of the $j$-th triangular factor $R_{j,j-1}$.
Obviously, in the implementation \cref{LE} is truncated to some large $T>0$.
In the end, each step of the DQR requires the solution of the IVPs \cref{IVPGAMMA} and the QR factorization \cref{GAMMAQR}.

\bigskip

The above summary was taken mainly from \cite{BredaDellaSchiava2018}, where the DQR is applied to the ODE obtained from the pseudospectral collocation of a given DDE (see \cref{sec:pscoll}), thus following the original approach of \cite{BredaDiekmannGyllenbergScarabelVermiglio2016} to also address the study of chaotic dynamics.
As anticipated in \cref{sec:intro}, the aim of the present work is to extend this procedure to more general classes of delay equations, such as REs and coupled equations.
Once the pseudospectral collocation is performed (possibly after linearization, see \cref{rem:lc-cl} later on), the outcome is an ODE like \cref{ODE}; thus, the DQR applies unchanged, independent of the original delay equation.

The literature on the theory and computation of LEs of ODEs is ample; for a starting reference, see \cite{DieciVanVleck2002}, but see also \cite{Adrianova1995} as a reference monograph.
QR methods were first proposed in the pioneering works \cite{BenettinGalganiGiorgilliStrelcyn1980a,BenettinGalganiGiorgilliStrelcyn1980b}; for a complete discussion of the discrete version, see \cite{DieciJollyVanVleck2011}.
The literature on the computation of LEs of delay equations is mostly of an experimental flavor and, to the best of the authors' knowledge, restricted only to DDEs.
As initial references, we can suggest \cite{BredaVanVleck2014,Breda2010,Farmer1982}, but see also \cite{ChekrounGhilLiuWang2016} for a more recent method to reduce DDEs to ODEs by using Galerkin-type projections.
Note also that all of these works rely on a Hilbert state space setting to legitimize orthogonal projections, while the technique in \cite{BredaDellaSchiava2018} is free from this constraint and thus maintains the classical state spaces (typically continuous functions for DDEs and absolutely integrable ones for REs).
Beyond the lack of relevant methods, this is part of the motivation for the extension of the approach proposed in \cite{BredaDellaSchiava2018} to more general delay equations.

\section{Pseudospectral collocation}
\label{sec:pscoll}

In this \lcnamecref{sec:pscoll}, we illustrate the use of pseudospectral collocation to reduce delay equations to ODEs, in view of the application of the DQR described in \cref{sec:dqr-ode}.
For the reader's convenience, we first present, separately, the discretization of an RE in \cref{sec:pscollX} and that of a DDE in \cref{sec:pscollY}, summarizing from, respectively, \cite{ScarabelDiekmannVermiglio2021} and \cite{BredaDellaSchiava2018} the main aspects for the present objective (for a full treatment, see again \cite{BredaDellaSchiava2018,ScarabelDiekmannVermiglio2021} and also \cite{BredaDiekmannGyllenbergScarabelVermiglio2016}).
Eventually, we combine the two approaches in \cref{sec:pscollXY} for a coupled equation.

In what follows, we use the subscripts $X$ and $Y$ to refer, respectively, to REs and DDEs.

\subsection{Pseudospectral collocation of REs}
\label{sec:pscollX}

Let $\tau>0$ be real and $d_{X}>0$ be an integer.
Consider the IVP for an RE given by
\begin{equation}\label{IVPX}
\left\{
\begin{alignedat}{2}
&x(t)=F(x_{t}),&&\quad t>0,\\
&x(\theta)=\phi(\theta),&&\quad\theta\in[-\tau,0],
\end{alignedat}
\right.
\end{equation}
where $\phi\in L^{1}\coloneqq L^{1}([-\tau,0];\mathbb{R}^{d_{X}})$, $F\colon L^{1}\to\mathbb{R}^{d_{X}}$ and $x_{t}$, defined as $x_{t}(\theta)\coloneqq x(t+\theta)$ for $\theta\in[-\tau,0]$, denotes the history or state function (so that $x_{0}=\phi$ represents the initial state).
If $F$ is globally Lipschitz, the IVP \cref{IVPX} has a unique solution on $[-\tau,+\infty)$ \cite[Theorem 3.8]{DiekmannGettoGyllenberg2008}.

In \cite{ScarabelDiekmannVermiglio2021} an efficient application of pseudospectral collocation to reduce \cref{IVPX} to an IVP for an ODE is proposed based on an equivalent formulation of \cref{IVPX} as an abstract Cauchy problem (ACP) describing the evolution of an integral of the original state $x_{t}$.
In particular, by defining the Volterra integral operator $\mathcal{V}\colon L^{1}\to AC_{0}$ as
\begin{equation*}
(\mathcal{V}\eta)(\theta)\coloneqq-\int_{\theta}^{0}\eta(s)\D s,
\end{equation*}
where $AC_{0}\coloneqq AC_{0}([-\tau,0];\mathbb{R}^{d_{X}})$ is the space%
\footnote{Hereinafter, we do not indicate the domain and codomain of a function space when clear from the context.}
of absolutely continuous functions vanishing at $0$, it turns out that \cref{IVPX} is equivalent to the ACP
\begin{equation}\label{ACPX}
\left\{
\begin{aligned}
&u'(t)=\mathcal{A}_{0,X}u(t)+q_{X}F(\mathcal{A}_{0,X}u(t)),\quad t\geq0,\\
&u(0)=\mathcal{V}\phi
\end{aligned}
\right.
\end{equation}
through $u(t)=\mathcal{V}x_{t}$.
Above, $\mathcal{A}_{0,X}\colon\mathcal{D}(\mathcal{A}_{0,X})\subset NBV_{0}\to NBV_{0}$ is defined as $(\mathcal{V}\restriction_{NBV_{0}})^{-1}$, i.e.,
\begin{equation}\label{A0X}
\mathcal{A}_{0,X}\mu\coloneqq\mu',\qquad \mathcal{D}(\mathcal{A}_{0,X})\coloneqq\{\mu\in AC_{0} : \mu=\mathcal{V}\eta\text{ for some }\eta\in NBV_{0}\},
\end{equation}
where $NBV_{0}$ is the space of functions of bounded variation vanishing at $0$ and continuous from the right, and $q_{X}\in NBV_{0}$ is defined as
\begin{equation*}
q_{X}(\theta)\coloneqq
\begin{cases}
0,&\theta=0,\\
-1,&\theta\in[-\tau,0).
\end{cases}
\end{equation*}

In order to discretize \cref{ACPX}, consider a mesh $\Omega_{M_{X},X}$ of $M_{X}$ points $-\tau\leq\theta_{M_{X},X}<\dots<\theta_{1,X}<0$ with $M_{X}$ a positive integer.
Correspondingly, let $P_{M_{X},X}\colon\mathbb{R}^{M_{X}d_{X}}\to NBV_{0}$ be the interpolation operator on $\{0\}\cup\Omega_{M_{X},X}$ with value $0$ at $\theta_{0,X}\coloneqq0$, i.e.,
\begin{equation*}
(P_{M_{X},X}\Phi)(\theta)\coloneqq\sum_{j=1}^{M_{X}}\ell_{j,X}(\theta)\Phi_{j},\quad\theta\in[-\tau,0],
\end{equation*}
where $\{\ell_{0,X},\ell_{1,X},\dots,\ell_{M_{X},X}\}$ is the Lagrange basis on $\{0\}\cup\Omega_{M_{X},X}$, and let $R_{M_{X},X}\colon NBV_{0}\to\mathbb{R}^{M_{X}d_{X}}$ be the restriction operator
\begin{equation*}
(R_{M_{X},X}\mu)_{j}\coloneqq\mu(\theta_{j,X}),\quad j=1,\dots,M_{X}.
\end{equation*}
Then, the discrete version of \cref{ACPX} is given by
\begin{equation}\label{ODEX}
\left\{
\begin{aligned}
&U'(t)=D_{M_{X},X}U(t)-\mathbf{1}_{M_{X},X}F_{M_{X}}(U(t)),\quad t\geq0,\\
&U(0)=R_{M_{X},X}\mathcal{V}\phi,
\end{aligned}
\right.
\end{equation}
with $U(t)\in\mathbb{R}^{M_{X}d_{X}}$ that approximates the integrated state $\mathcal{V}x_{t}$ according to $U_{j}(t)\approx(\mathcal{V}x_{t})(\theta_{j,X})$, $j=1,\dots,M_{X}$%
\footnote{Note then that $\mathcal{A}_{0,X}P_{M_{X},X}U(t)\approx x_{t}$.}%
, and where $D_{M_{X},X}\coloneqq R_{M_{X},X}\mathcal{A}_{0,X}P_{M_{X},X}\in\mathbb{R}^{M_{X}d_{X}\times M_{X}d_{X}}$ has $d_{X}\times d_{X}$-block entries
\begin{equation*}
[D_{M_{X},X}]_{i,j}=\ell'_{j,X}(\theta_{i,X})I_{d_{X}},\quad i,j=1,\dots,M_{X},
\end{equation*}
where $I_{d_{X}}$ is the identity on $\mathbb{R}^{d_{X}}$, $F_{M_{X}}\coloneqq F\circ\mathcal{A}_{0,X}P_{M_{X},X}$ and $\mathbf{1}_{M_{X},X}\in\mathbb{R}^{M_{X}d_{X} \times d_{X}}$ has all $d_{X}\times d_{X}$-block entries $I_{d_{X}}$%
\footnote{Note that $\mathbf{1}_{M_{X},X}$ discretizes $-q_{X}$.}%
.

\begin{remark}\label{rem:sunstarX}
Instead of $NBV_{0}$, \cite{ScarabelDiekmannVermiglio2021} uses the space $NBV$ of functions of bounded variation that vanish at $0$ and are continuous from the right on $(-\tau, 0)$, but not necessarily at $-\tau$.
In that setting, $\mathcal{C}_{0,X}\coloneqq(\mathcal{V}\restriction_{NBV})^{-1}$ is a multi-valued operator, defined as%
\footnote{It turns out that $D(\mathcal{C}_{0,X})=D(\mathcal{A}_{0,X})$ (see \cref{A0X}).}
$\mathcal{C}_{0,X}\mu \coloneqq \{\eta:\mu=V\eta\}$, since functions differing only by the jump at $-\tau$ are mapped by $\mathcal{V}$ to the same element of $NBV$.
The \emph{trivial} semigroup $\{S_{0,X}(t)\}_{t\geq0}$ on $NBV$ defined as
\begin{equation*}
S_{0,X}(t)\colon NBV\to NBV,
\qquad
(S_{0,X}(t)\mu)(\theta)\coloneqq
\begin{cases}
\mu(t+\theta),&t+\theta\leq0,\\
0,&t+\theta>0,
\end{cases}
\end{equation*}
is not strongly continuous.
However, its restriction $\{T_{0,X}(t)\}_{t\geq0}$ to $AC_{0}$ is strongly continuous and $\mathcal{A}_{0,X}$ is its infinitesimal generator.

From this point of view, the semilinear ACP \cref{ACPX} renders a clear separation between the \emph{translation along the solutions} (through the linear semigroup $\{T_{0,X}(t)\}_{t\geq0}$ and its infinitesimal generator $\mathcal{A}_{0,X}$) and the \emph{rule for extension} (basically through the nonlinear right-hand side $F$ of the specific RE), which are the two ingredients of a delay equation.
For these and related aspects of the theory of delay equations, see \cite{DiekmannGettoGyllenberg2008,DiekmannVanGilsVerduynLunelWalther1995} for the sun--star ($\odot\ast$) theory and \cite{DiekmannVerduynLunel2021,ScarabelDiekmannVermiglio2021} for the more recent twin semigroup theory.
\end{remark}

\subsection{Pseudospectral collocation of DDEs}
\label{sec:pscollY}

Let $\tau>0$ be real and $d_{Y}>0$ be an integer.
Consider the IVP for a DDE given by
\begin{equation}\label{IVPY}
\left\{
\begin{alignedat}{2}
&y'(t)=G(y_{t}),&&\quad t\geq0,\\
&y(\theta)=\psi(\theta),&&\quad\theta\in[-\tau,0],
\end{alignedat}
\right.
\end{equation}
where $\psi\in C\coloneqq C([-\tau,0];\mathbb{R}^{d_{Y}})$, $G\colon C\to\mathbb{R}^{d_{Y}}$ and $y_{t}$ is defined as $x_{t}$ in \cref{sec:pscollX} (so, again, $y_{0}=\psi$ represents the initial state).
If $G$ is globally Lipschitz, the IVP \cref{IVPY} has a unique solution on $[-\tau,+\infty)$ \cite[Section 2.2]{HaleVerduynLunel1993}.

In \cite{BredaDellaSchiava2018} (but see also \cite{BredaDiekmannGyllenbergScarabelVermiglio2016}) pseudospectral collocation is used to reduce \cref{IVPY} to an IVP for an ODE based on an equivalent formulation of \cref{IVPY} as an ACP describing the evolution of the original state $y_{t}$, viz.
\begin{equation}\label{ACPY}
\left\{
\begin{aligned}
&v'(t)=\mathcal{A}_{Y}v(t),\quad t\geq0,\\
&v(0)=\psi
\end{aligned}
\right.
\end{equation}
through $v(t)=y_{t}$.
Above, $\mathcal{A}_{Y}\colon\mathcal{D}(\mathcal{A}_{Y})\subset C\to C$ is defined as
\begin{equation*}
\mathcal{A}_{Y}\rho\coloneqq\rho',\qquad \mathcal{D}(\mathcal{A}_{Y})\coloneqq\{\rho\in C : \rho'\in C\text{ and }\rho'(0)=G(\rho)\}.
\end{equation*}

In order to discretize \cref{ACPY}, consider a mesh $\Omega_{M_{Y},Y}$ of $1+M_{Y}$ points $-\tau=\theta_{M_{Y},Y}<\theta_{M_{Y}-1,Y}<\dots<\theta_{1,Y}<\theta_{0,Y}\coloneqq0$ with $M_{Y}$ a positive integer.
Correspondingly, let $P_{M_{Y},Y}\colon\mathbb{R}^{(1+M_{Y})d_{Y}}\to C$ be the interpolation operator on $\Omega_{M_{Y},Y}$, i.e.,
\begin{equation*}
(P_{M_{Y},Y}\Psi)(\theta)\coloneqq\sum_{j=0}^{M_{Y}}\ell_{j,Y}(\theta)\Psi_{j},\quad\theta\in[-\tau,0],
\end{equation*}
where $\{\ell_{0,Y},\ell_{1,Y},\dots,\ell_{M_{Y},Y}\}$ is the Lagrange basis on $\Omega_{M_{Y},Y}$, and let $R_{M_{Y},Y}:C\to\mathbb{R}^{(1+M_{Y})d_{Y}}$ be the restriction operator
\begin{equation*}
(R_{M_{Y},Y}\rho)_{j}\coloneqq\rho(\theta_{j,Y}),\quad j=0,1,\dots,M_{Y}.
\end{equation*}
Then, the discrete version of \cref{ACPY} is given by
\begin{equation}\label{ODEY}
\left\{
\begin{alignedat}{2}
&V_{0}'(t)=G_{M_{Y}}(V(t)),&&\quad t\geq0,\\
&V_{j}'(t)=D_{M_{Y},Y}V(t),&&\quad t\geq0,\; j=1,\dots,M_{Y},\\
&V(0)=R_{M_{Y},Y}\psi,
\end{alignedat}
\right.
\end{equation}
for $V(t)\in\mathbb{R}^{(1+M_{Y})d_{Y}}$ that approximates the state $y_{t}$ according to $V_{j}(t)\approx(y_{t})(\theta_{j,Y})$, $j=0,1,\dots,M_{Y}$, and where $D_{M_{Y},Y}$ has $d_{Y}\times d_{Y}$-block entries
\begin{equation*}
[D_{M_{Y},Y}]_{i,j}=\ell'_{j,Y}(\theta_{i,Y})I_{d_{Y}},\quad i=1,\dots,M_{Y},\;j=0,1,\dots,M_{Y},
\end{equation*}
where $I_{d_{Y}}$ is the identity on $\mathbb{R}^{d_{Y}}$ and $G_{M_{Y}}\coloneqq G\circ P_{M_{Y},Y}$.

\begin{remark}\label{rem:sunstarY}
Note that the ACP \cref{ACPY} can be alternatively described as the equivalent semilinear ACP
\begin{equation}\label{ACPYss}
\left\{
\begin{aligned}
&v'(t)=\mathcal{A}_{0,Y}v(t)+q_{Y}G(v(t)),\quad t\geq0,\\
&v(0)=\psi,
\end{aligned}
\right.
\end{equation}
where $\mathcal{A}_{0,Y}$ is the infinitesimal generator of the \emph{shift} semigroup $\{T_{0,Y}(t)\}_{t\geq0}$ defined as
\begin{equation*}
T_{0,Y}(t)\colon C\to C,
\qquad
(T_{0,Y}(t)\rho)(\theta)\coloneqq
\begin{cases}
\rho(t+\theta),&t+\theta\leq0,\\
\rho(0),&t+\theta>0,
\end{cases}
\end{equation*}
and $q_{Y}\in L^{\infty}$ is defined as
\begin{equation*}
q_{y}(\theta)\coloneqq
\begin{cases}
1,&\theta=0,\\
0,&\theta\in[-\tau,0).
\end{cases}
\end{equation*}
\Cref{ACPYss} renders for DDEs the same separation between \emph{translation along the solutions} and \emph{rule for extension} as illustrated in \cref{rem:sunstarX} for REs (see again \cite{DiekmannVanGilsVerduynLunelWalther1995}).
The pseudospectral collocation of \cref{ACPYss} leads, again, to \cref{ODEY}, which can be rewritten equivalently as
\begin{equation}\label{ODEYss}
\left\{
\begin{aligned}
&V'(t)=D_{0,M_{Y},Y}V(t)+\mathbf{1}_{M_{Y},Y}G_{M_{Y}}(U(t),V(t)),\quad t\geq0,\\
&V(0)=R_{M_{Y},Y}\psi,
\end{aligned}
\right.
\end{equation}
where $D_{0,M_{Y},Y}$ is as $D_{M_{Y},Y}$ but with an additional $d_{Y}$-block row of zeros; also, $\mathbf{1}_{M_{Y},Y}\in\mathbb{R}^{(1+M_{Y})d_{Y}\times d_{Y}}$ has the first $d_{Y}\times d_{Y}$-block equal to $I_{d_{Y}}$ and all of the others equal to zero.
Now, \cref{ODEYss} resembles \cref{ACPYss}.
\end{remark}

\subsection{Pseudospectral collocation of coupled equations}
\label{sec:pscollXY}

Let $\tau$, $d_{X}$ and $d_{Y}$ be as above.
Consider the IVP for a coupled equation given by
\begin{equation}\label{IVPXY}
\left\{
\begin{alignedat}{2}
&x(t)=F(x_{t},y_{t}),&&\quad t>0,\\
&y'(t)=G(x_{t},y_{t}),&&\quad t\geq0,\\
&x(\theta)=\phi(\theta),&&\quad\theta\in[-\tau,0],\\
&y(\theta)=\psi(\theta),&&\quad\theta\in[-\tau,0],
\end{alignedat}
\right.
\end{equation}
where $\phi\in L_{1}$, $\psi\in C$, $F\colon L^{1}\times C\to\mathbb{R}^{d_{X}}$ and $G\colon L^{1}\times C\to\mathbb{R}^{d_{Y}}$.
For well-posedness, see \cite{DiekmannGettoGyllenberg2008}.

By combining the approaches of the previous \lcnamecrefs{sec:pscollX}, it follows that \cref{IVPXY} is equivalent to the ACP
\begin{equation}\label{ACPXY}
\left\{
\begin{aligned}
&(u'(t),v'(t))=\mathcal{B}_{X,Y}(u(t),v(t)),\quad t\geq0,\\
&(u(0),v(0))=(\mathcal{V}\phi,\psi),
\end{aligned}
\right.
\end{equation}
with $\mathcal{B}_{X,Y}\colon\mathcal{D}(\mathcal{B}_{X,Y})\subset NBV_{0}\times Y\to NBV_{0}\times Y$ defined as
\begin{gather*}
\mathcal{B}_{X,Y}(\phi,\psi)\coloneqq(\mathcal{A}_{0,X}\phi+q_{X}F(\mathcal{A}_{0,X}\phi,\psi),\psi'),\\
\mathcal{D}(\mathcal{B}_{X,Y})\coloneqq\{(\varphi,\psi)\in\mathcal{D}(\mathcal{A}_{0,X})\times Y : \psi'\in Y,\ \psi'(0)=G(\mathcal{A}_{0,X}\phi,\psi)\},
\end{gather*}
through $(u(t),v(t))=(\mathcal{V}x_{t},y_{t})$.

The discrete version of \cref{ACPXY} is as follows:
\begin{equation}\label{ODEXY}
\left\{
\begin{alignedat}{2}
&U'(t)=D_{M_{X},X}U(t)-\mathbf{1}_{M_{X},X}F_{M_{X}}(U(t),V(t)),&\quad& t\geq0,\\
&V_{0}'(t)=G_{M_{Y}}(U(t),V(t)),&\quad& t\geq0,\\
&V_{j}'(t)=D_{M_{Y},Y}V(t),&\quad& t\geq0,\; j=1,\dots,M_{Y},\\
&U(0)=R_{M_{X},X}\mathcal{V}\phi,\\
&V(0)=R_{M_{Y},Y}\psi.
\end{alignedat}
\right.
\end{equation}
Note that, now, $F_{M_{X}}\coloneqq F\circ(\mathcal{A}_{0,X}P_{M_{X},X},P_{M_{Y},Y})$ and $G_{M_{Y}}\coloneqq G\circ(\mathcal{A}_{0,X}P_{M_{X},X},P_{M_{Y},Y})$.
The total number of approximating ODEs is $M_{X}d_{X}+(1+M_{Y})d_{Y}$, which becomes $(2M+1)d$ if $d_{X}=d_{Y}=d$ and $M_{X}=M_{Y}=M$.

\begin{remark}
Following \cref{rem:sunstarX,rem:sunstarY}, it is not difficult to see that \cref{ACPXY} is equivalent to the semilinear ACP
\begin{equation}\label{ACPXYss}
\left\{
\begin{aligned}
&(u'(t),v'(t))=\mathcal{A}_{0,X,Y}(u(t),v(t))+\mathcal{N}_{X,Y}(u(t),v(t)),\quad t\geq0,\\
&(u(0),v(0))=(\mathcal{V}\phi,\psi),
\end{aligned}
\right.
\end{equation}
where $\mathcal{A}_{0,X,Y}\coloneqq\diag(\mathcal{A}_{0,X},\mathcal{A}_{0,Y})$ is linear and
\begin{equation*}
\mathcal{N}_{X,Y}(\phi,\psi)\coloneqq(q_{X}F(\mathcal{A}_{0,X}\phi,\psi),q_{Y}G(\mathcal{A}_{0,X}\phi,\psi))
\end{equation*}
is nonlinear.
The pseudospectral collocation of \cref{ACPXYss} leads, again, to \cref{ODEXY}, where, correspondingly, the ODEs for $V$ can be compacted as done in \cref{rem:sunstarY} for DDEs; see \cref{ODEYss}.
The (numerical) analysis of \cref{ACPXYss} is current work in progress at CDLab%
\footnote{\url{http://cdlab.uniud.it/}}%
,
also in view of the corresponding sun--star theory of coupled equations developed in \cite{DiekmannGettoGyllenberg2008} and of the more recent twin semigroup theory of \cite{DiekmannVerduynLunel2021}.
\end{remark}

\section{Implementation}
\label{sec:impl}

Due to our choice of example equations (see \cref{sec:res}), in our implementation we considered scalar equations ($d_X=d_Y=1$) of the following kinds:
\begin{equation}\label{REprototype}
x(t) = \int_{-\tau_2}^{-\tau_1} f(x(t+\theta)) \D\theta
\end{equation}
for REs, and
\begin{equation}\label{coupledprototype}
\left\{
\begin{aligned}
&x(t) = y(t) \int_{-\tau_2}^{-\tau_1} f_1(x(t+\theta)) \D\theta, \\
&y'(t) = g(y(t)) + y(t) \int_{-\tau_2}^{-\tau_1} f_2(x(t+\theta)) \D\theta,
\end{aligned}
\right.
\end{equation}
for coupled equations%
\footnote{Observe that the unknown of the differential equation is not delayed.},
where $\tau_2>\tau_1>0$, and $f$, $f_1$, $f_2$ and $g$ are (possibly) nonlinear functions $\mathbb{R}\to\mathbb{R}$.
Nevertheless, generalizing to other forms of equations is usually fairly straightforward%
\footnote{%
For instance, one can consider models such as (3.3) in \cite{RipollFont2023}, yet with finite age-span, which do not enter class \cref{coupledprototype}.
In such cases, the method is implemented following the same strategy described here (discretization of the nonlinear system for computing solutions, linearization and discretization of the linearized system for computing the LEs).
However, the authors are currently implementing a general code for a larger class of equations (ideally the most general \cref{IVPXY}), which was beyond the scope of this work.
}%
.

We implemented the pseudospectral discretization%
\footnote{For more details on pseudospectral methods, see also \cite{Trefethen2000}.}
using Chebyshev nodes of type II (extrema) as the meshes $\{0\}\cup\Omega_{M_X,X}$ and $\Omega_{M_Y,Y}$ of points in $[-\tau, 0]$ with $M_X=M_Y$.
To compute the nodes and the corresponding differentiation matrix, we used the \texttt{cheb} routine of \cite[Chapter 6]{Trefethen2000}.
For the interpolation, we used the barycentric Lagrange interpolation formula \cite{BerrutTrefethen2004}; the barycentric weights corresponding to Chebyshev extrema are explicitly known and are given therein.
For the quadrature of the integrals, we used the Clenshaw--Curtis formula \cite{ClenshawCurtis1960,Trefethen2008}.
We implemented $(\mathcal{A}_{0,X}P_{M_{X},X}\Phi)(\theta)$ as $\sum_{j=1}^{M_{X}}\ell_{j,X}'(\theta)\Phi_{j}$, computing $\ell_{j,X}'$ as the polynomial interpolating the $j$-th column of the differentiation matrix, again with the barycentric formula.

In order to apply the DQR described in \cref{sec:dqr-ode} to the approximating ODE, the latter needs to be linearized around a reference solution.
The linearization is done explicitly.
The solutions are computed by using MATLAB's \texttt{ode45}, which implements the embedded Dormand--Prince $(5,4)$ method \cite{DormandPrince1980,ShampineReichelt1997}.
For the differential part of \cref{coupledprototype}, the initial value consists of the vector of values of the chosen initial function at $\Omega_{M_Y,Y}$, while, for \cref{REprototype} and the renewal part of \cref{coupledprototype}, the vector $D_{M_X,X}^{-1}u$ is used%
\footnote{In \cref{ODEX} the initial value is specified as $R_{M_{X},X}\mathcal{V}\phi$, i.e., the vector representing the polynomial interpolating the exact integral of the initial value $\psi$.
Another approach is to use $R_{M_{X},X}\mathcal{V}P_{M_{X},X}R_{M_{X},X}\phi$, in which the integral of the polynomial interpolating $\psi$ is used.
In our implementation we use neither; our choice is computationally easier and is motivated as follows.

As already noted, in order to represent the integrated state, only the vector $U$ of values at $\Omega_{M_X,X}$ is needed, as the value at $\theta_{0,X}=0$ is always $0$.
Computing the derivative of the interpolating polynomial by applying the differentiation matrix to $(0, U)^T$ (where the $0$ stands for a column vector of $d_X$ zeros), we obtain $(d_{M_X,X} U, D_{M_X,X} U)^T$, where $d_{M_X,X}\in\mathbb{R}^{d_{X}\times M_{X}d_{X}}$ is a row of $d_{X}\times d_{X}$-block entries $\ell'_{j,X}(\theta_{0,X})I_{d_{X}}$ for $j=1,\dots,M_{X}$.
Since deriving a polynomial lowers its degree by one, $D_{M_X,X} U$ uniquely determines the derivative of the polynomial represented by $U$, which motivates our use of $D_{M_X,X}^{-1}u$.}%
, where $u$ is the vector of values of the initial function at $\Omega_{M_X,X}$.

Finally, the DQR for a linear ODE is implemented in the \texttt{dqr} routine of \cite{BredaDellaSchiava2018}, which follows \cite{DieciVanVleck2004}.
Therein, the IVPs \cref{IVPGAMMA} are again solved with the Dormand--Prince $(5,4)$ pair; however, instead of adapting the step size (initially $0.01$) based on the error between the two solutions, the automatic adaptation controls the error between the corresponding LEs.
As an initial guess for the fundamental matrix solution, a random matrix is used.
The computation is stopped when the specified truncation time $T$ is reached.

\begin{remark}\label{rem:lc-cl}
For REs only, and in particular for the example described in \cref{sec:quadRE}, we experimented also with a different method, based on computing a solution of \cref{REprototype}, linearizing the latter around the former and applying the pseudospectral collocation to the resulting linear RE.
We computed the solution of the RE with the method described in \cite{MessinaRussoVecchio2008}, which is based on the trapezoidal quadrature formula on a uniform grid in $[-\tau_2, 0]$ with the constraint that $-\tau_1$ must be a grid point.
Corresponding to a solution $\bar{x}$, we considered the linear RE%
\footnote{\label{fn:linRE}
In general, \cref{REprototypelin} may actually not be the linearization of \cref{REprototype} around $\bar{x}$ in $L^1$.
Indeed, the right-hand side of the equation is \emph{not} Fr\'echet-differentiable unless $f$ is affine.
See \cite[Section 3.5]{DiekmannGettoGyllenberg2008} for details, in particular with respect to studying the stability of equilibria; the extension of the results therein is an open problem.}
\begin{equation}\label{REprototypelin}
x(t) = \int_{-\tau_2}^{-\tau_1} f'(\bar{x}(t+\theta)) x(t+\theta) \D\theta.
\end{equation}
See \cref{sec:quadRE} for a comparison of the approaches.
\end{remark}

\begin{remark}\label{rem:linRE}
In \cref{rem:lc-cl} (more precisely in \cref{fn:linRE}), we observed that, in most cases, REs cannot be linearized.
However, in many of those cases, the ODE resulting from the pseudospectral discretization can, in fact, be linearized; for example, the ODEs resulting from \cref{REprototype} and \cref{coupledprototype} can be linearized if $f$, $f_1$, $f_2$ and $g$ are differentiable.
As an example, the linearization of \cref{ODEXY} around a solution $(\overline{U}, \overline{V})$ is as follows:
\begin{equation*}
\left\{
\begin{aligned}
&U'(t)=D_{M_{X},X}U(t)-\mathbf{1}_{M_{X},X} \cdot JF_{M_{X}}(\overline{U}(t),\overline{V}(t)) \cdot (U(t),V(t))^{T},\\
&V_{0}'(t)=JG_{M_{Y}}(\overline{U}(t),\overline{V}(t)) \cdot (U(t),V(t))^{T},\\
&V_{j}'(t)=D_{M_{Y},Y}V(t),\quad j=1,\dots,M_{Y},
\end{aligned}
\right.
\end{equation*}
where $J$ indicates the Jacobian matrix, $JF_{M_{X}}(\overline{U}(t),\overline{V}(t))$ is a $d_X \times (M_Xd_X+(1+M_Y)d_Y)$ matrix and $JG_{M_{Y}}(\overline{U}(t),\overline{V}(t))$ is a $d_Y \times (M_Xd_X+(1+M_Y)d_Y)$ matrix.
In \cref{sec:quadRE} below, we explicitly show the linearized ODE for an example RE.
\end{remark}

We recall that the MATLAB codes implementing the method and the scripts to reproduce the experiments of \cref{sec:res} are available at \url{http://cdlab.uniud.it/software}.

\section{Results}
\label{sec:res}

We present here three example equations: an RE with a quadratic (logistic-like) nonlinearity in \cref{sec:quadRE}, an RE modeling egg cannibalism in \cref{sec:cannibRE} and a simplified version of the Daphnia model with a logistic term for the growth of the resource in \cref{sec:logDaphnia}.
In particular, we use the first example to test the proposed method also from the numerical point of view; we then apply it to the second and third example to compute the exponents.

\subsection{RE with quadratic nonlinearity}
\label{sec:quadRE}

The first equation we study is the RE with quadratic nonlinearity from \cite{BredaDiekmannLiessiScarabel2016}:
\begin{equation}\label{quadRE}
x(t) = \frac{\gamma}{2}\int_{-3}^{-1} x(t+\theta)(1-x(t+\theta)) \D\theta,
\end{equation}
i.e., \cref{REprototype} with $\tau_1\coloneqq1$, $\tau_2\coloneqq3$ and $f(x) \coloneqq \frac{\gamma}{2}x(1-x)$.
Its equilibria and their stability properties are known; in particular, its nontrivial equilibrium undergoes a Hopf bifurcation for $\gamma = 2 + \frac{\pi}{2}$ and the branch of periodic solutions arising from there has the analytic expression
\begin{equation}\label{quadREpersol}
\bar{x}(t) = \frac{1}{2} + \frac{\pi}{4 \gamma} + \sqrt{\frac{1}{2} - \frac{1}{\gamma} - \frac{\pi}{2 \gamma^{2}}\Bigl(1+\frac{\pi}{4}\Bigr)} \sin\Bigl(\frac{\pi}{2} t\Bigr).
\end{equation}
Observe that the period is $4$, independent of $\gamma$.
Moreover, it is experimentally known that it presents several period-doubling bifurcations, possibly leading to a cascade and, eventually, to chaos \cite{BredaDiekmannLiessiScarabel2016}.

Since several properties of \cref{quadRE} are analytically known, we use it to test the effectiveness and efficiency of the method proposed for LE computation, and to compare it to the alternative approach described in \cref{rem:lc-cl}.

For equilibria, the LEs are the real parts of the eigenvalues $\lambda$ of the infinitesimal generator of the semigroup of solution operators, which are related to the eigenvalues $\mu$ of the solution operator that advances the solution by $h$ via
\begin{equation}\label{mulambda}
\mu = \E^{\lambda h}.
\end{equation}
For periodic solutions, the LEs are the real parts of the Floquet exponents, which are related to the Floquet multipliers (i.e., the eigenvalues of the monodromy operator) via \cref{mulambda}, where $\mu$, $\lambda$ and $h$ are, respectively, a Floquet multiplier, a Floquet exponent and the period.
In both cases, we can thus obtain the LEs by computing the eigenvalues $\mu$ of an evolution operator with any time step $h$ for the equilibria (we choose $h=\tau_2=3$ for \cref{quadRE}) and a time step $h$ equal to the period for periodic solutions ($h=4$ for \cref{quadREpersol}), and then computing the real part of $\log(\mu)/h$.
In order to obtain reference values for our experiments, we compute the spectra of evolution operators with the method of \cite{BredaLiessi2018}, which is based on the pseudospectral collocation of the operator; we use the implementation \texttt{eigTMNpw} of \cite{BredaLiessiVermiglio2022,BredaLiessiVermiglio}.

\bigskip

Although computing the solutions of delay equations is not the focus of this work, given that both the main approach and the alternative one of \cref{rem:lc-cl} involve computing solutions, our first experiment compares the error of the computed solutions with respect to the known periodic solutions of \cref{quadRE}.
We choose $\gamma = 4>2+\frac{\pi}{2}$, which corresponds to a stable periodic solution, since the first period-doubling bifurcation is experimentally known to happen at $\gamma \approx 4.32$ \cite{BredaDiekmannLiessiScarabel2016,BredaLiessi2018}.

\begin{figure}
\begin{center}
\includegraphics[trim=24.254pt 0 88.885pt 0]{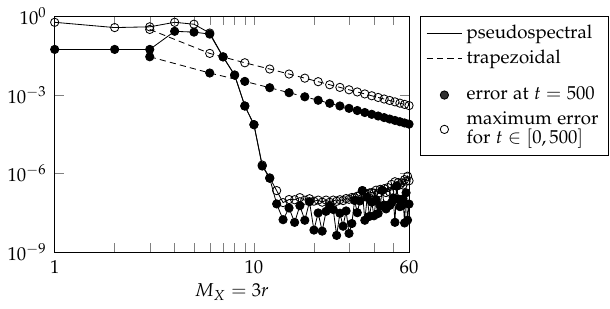}
\caption{Errors on the solution of the RE with quadratic nonlinearity \cref{quadRE} with $\gamma=4$ with respect to the known periodic solution \cref{quadREpersol}, computed via pseudospectral discretization (solid lines) and directly with the trapezoidal method (dashed lines), measured as the absolute error at $t=500$ ($\bullet$) and as the maximum absolute error on a grid of points in $[0, 500]$ ($\circ$), when varying the number of nodes (minus $1$) in the grid in $[-3, 0]$, i.e., $M_X$ for the pseudospectral discretization and $3r$ for the trapezoidal method.
The exact periodic solution \cref{quadREpersol} is used as the initial value.}
\label{fig:quad-sol-compare}
\end{center}
\end{figure}

\Cref{fig:quad-sol-compare} shows the errors on the solution of the approximating ODE and on the solution of the original RE \cref{quadRE} with respect to the number of nodes (minus $1$) in the grid in $[-3, 0]$, i.e., $M_X$ for the pseudospectral discretization and $3r$ for the trapezoidal method%
\footnote{%
As noted in \cref{rem:lc-cl}, $-\tau_1=-1$ must be a grid point in $[-\tau_2,0]=[-3,0]$; we thus choose the number $r$ of nodes (minus $1$) in $[-\tau_1,0]$ as the discretization parameter for the trapezoidal method \cite{MessinaRussoVecchio2008}, resulting in $3r+1$ nodes in total.
}
\cite{MessinaRussoVecchio2008}; in both cases, two errors are measured, namely the absolute error at $t = 500$ and the maximum absolute error on a grid of points in $[0, 500]$ (a uniform grid with step $0.05$ for the pseudospectral approach, the time points given by the trapezoidal method for the alternative approach).
To solve the approximating ODE given by the pseudospectral discretization, we used \texttt{ode45} with $\texttt{RelTol} = 10^{-6}$ and $\texttt{AbsTol} = 10^{-7}$, which justifies the barrier on the error in \cref{fig:quad-sol-compare}.

The experiment confirms that the trapezoidal method has order $2$, as proved in \cite{MessinaRussoVecchio2008}, and that the pseudospectral discretization has infinite order, which is often the case for pseudospectral methods applied to smooth problems \cite{Trefethen2000}.
Even for rather small values of $M_X=3r$, the error for the pseudospectral method is several orders of magnitudes smaller than the one for the trapezoidal method%
\footnote{As already noted, computing the solutions is not the focus of the present work.
Admittedly, there are other more sophisticated methods in the literature: see, e.g., \cite{BellenJackiewiczVermiglioZennaro1989,Vermiglio1992,Brunner1994b}.
However, in most cases, they are not readily applicable to \cref{REprototype,coupledprototype}, due to the discontinuity in the integration kernel at $-\tau_1$, when the integral is considered on the interval $[-\tau_2,0]$; in other cases, the implementation of the method is not available and is not as straightforward as \cite{MessinaRussoVecchio2008}.
It is worth mentioning that pseudospectral methods for computing periodic solutions of REs and coupled equations, exhibiting the usual infinite order of convergence, are available in \cite{Ando2021,AndoBreda2023,AndoBreda}; they are based on solving the corresponding boundary value problem, so they cannot be straightforwardly adapted to computing generic solutions.}%
.

\bigskip

In the next experiment, we investigate how the errors on the LEs depend on the choice of $M_X$ and of the final time%
\footnote{Observe that the \texttt{dqr} routine does not stop exactly at $T$, but at the first time step exceeding $T$, i.e., it does not refine the final step.}
$T$.
We choose values of $\gamma$ corresponding to the stable trivial equilibrium ($\gamma=0.5$), the stable nontrivial equilibrium ($\gamma=3$) and the stable periodic orbit ($\gamma=4$).

Since we are going to use the linearization of the ODE \cref{ODEX} coming from the RE \cref{quadRE}, as an example, we show it explicitly here.
With reference to \cref{rem:linRE}, observe that the right-hand side of \cref{quadRE} is not Fr\'echet-differentiable as a map from $L^1$ to $\mathbb{R}$, while the right-hand side of the discretized equation is differentiable.
The linearization of the approximating ODE around the solution $\overline{U}$ is given by
\begin{equation*}
U'(t)=D_{M_{X},X}U(t)-\mathbf{1}_{M_{X},X} \cdot JF_{M_{X}}(\overline{U}(t)) \cdot U(t),
\end{equation*}
where $JF_{M_{X}}(\overline{U}(t))$ is a row vector with components
\begin{equation*}
[JF_{M_{X}}(\overline{U}(t))]_{j} = \frac{\gamma}{2}\int_{-3}^{-1} \biggl(1-2\sum_{k=1}^{M_{X}}\ell_{k,X}'(\theta)\overline{U}_{k}(t)\biggr) \ell_{j,X}'(\theta) \D\theta, \quad j=1,\dots,M_X.
\end{equation*}

\begin{figure}
\begin{center}
\includegraphics[trim=24.258pt 0 91.917pt 0]{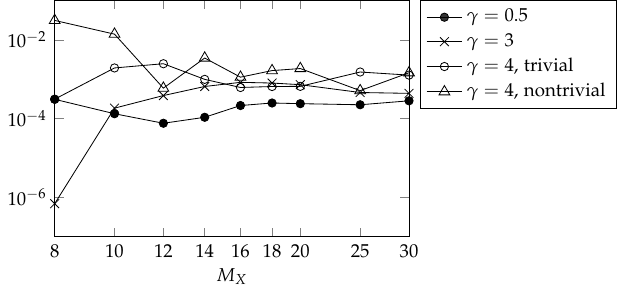}
\caption{Absolute errors on the dominant LEs of the RE with quadratic nonlinearity \cref{quadRE} for values of $\gamma$ corresponding to the stable trivial equilibrium ($\gamma=0.5$), the stable nontrivial equilibrium ($\gamma=3$) and the stable periodic orbit ($\gamma=4$).
For the last one, both the trivial and the dominant nontrivial exponents are shown.
The errors are measured with respect to the exponents computed via \texttt{eigTMNpw}.
The final time for \texttt{dqr} is $T=1000$.}
\label{fig:quad-le-M}
\end{center}
\end{figure}

\begin{figure}
\begin{center}
\includegraphics[trim=24.254pt 0 10.026pt 0]{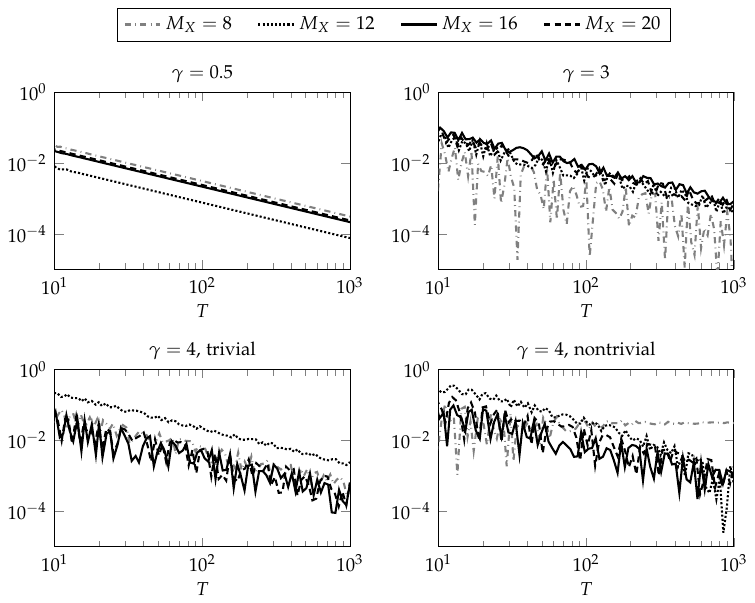}
\caption{Absolute errors on the dominant LEs of the RE with quadratic nonlinearity \cref{quadRE} for values of $\gamma$ corresponding to the stable trivial equilibrium ($\gamma=0.5$), the stable nontrivial equilibrium ($\gamma=3$) and the stable periodic orbit ($\gamma=4$).
For the last one, both the trivial and the dominant nontrivial exponents are shown.
The errors are measured with respect to the exponents computed via \texttt{eigTMNpw}.
The exponents are computed for $M_X \in \{8, 12, 16, 20\}$ as shown in the legend.}
\label{fig:quad-le-T-1000}
\end{center}
\end{figure}

In \cref{fig:quad-le-M,fig:quad-le-T-1000}, we can see the absolute errors on the dominant LE increasing either $M_X$ or $T$.
The tolerance for \texttt{dqr} is $10^{-6}$, while those for \texttt{ode45} are $\texttt{RelTol} = 10^{-6}$ and $\texttt{AbsTol} = 10^{-7}$.
The reference values are obtained by using \texttt{eigTMNpw} with the default options and $120$ as the degree of the collocation polynomials (fixed independently of $M_X$).

In \cref{fig:quad-le-M} the final time $T=1000$ is fixed and $M_X$ increases.
We can observe that, apart from very low values of $M_X$, the error reaches a barrier%
\footnote{%
Note in \cref{fig:quad-le-M} that for $\gamma=3$ the error is initially below the final barrier; we do not have an explanation of this phenomenon, but in \cref{fig:quad-le-T-1000} for $\gamma=3$ and $M_X=8$ we can observe rather erratic behavior of the error when varying $T$, despite the fact that it has a linearly decreasing bound.
Moreover, we remark that the value of $M_X$ required to reach the error barrier for a given $T$ depends, in general, on the specific equation and the specific values of its parameters.
}%
.
We performed the same experiment with $T=10000$ and could make the same observation, although the barrier was smaller by about one order of magnitude.
The barrier depends on the error due to the time truncation in \cref{LE}.
Indeed, \cref{fig:quad-le-T-1000}, where $M_X$ is constant and $T$ varies, shows that the LEs converge linearly (confirming what is explained in \cite{BredaVanVleck2014}).
In \cref{fig:quad-le-M} the truncation error appears to dominate on the error due to the collocation.

For the dominant nontrivial exponent%
\footnote{The trivial exponent $0$ is always an LE for a periodic solution due to the translation invariance.}
of the periodic solution, we observe in \cref{fig:quad-le-T-1000} that, for $M_X=8$, the error seems to reach a barrier, indicating that more ODEs are necessary to reproduce the properties of the original RE more accurately, as it is reasonable to expect.
In other experiments, we observed that, as $T$ increases, error barriers are reached also for increasing values of $M_X$.

We showed here the results of the main approach only.
We performed the same experiments also with the alternative approach of \cref{rem:lc-cl}, using, for the linearization, both the exact solutions (which are known for the chosen values of $\gamma$) and the numerical solutions computed by using the trapezoidal method.
With the exact solutions, we obtained almost exactly the same values: this means that the solution of the ODE is a good enough approximation of the solution of the RE, and that exchanging the linearization and the collocation does not influence the results.
However, when using the numerical solutions obtained via the trapezoidal method, the errors on the LEs were higher: in the experiment shown in \cref{fig:quad-le-M}, the errors in the periodic case were one order of magnitude larger, while in the experiment of \cref{fig:quad-le-T-1000} the errors reached barriers ranging between $10^{-2}$ and $10^{-1}$.
For these reasons, we henceforth use only the main approach as the more practical one.

\begin{figure}
\begin{center}
\includegraphics[trim=24.937pt 0 9.409pt 0]{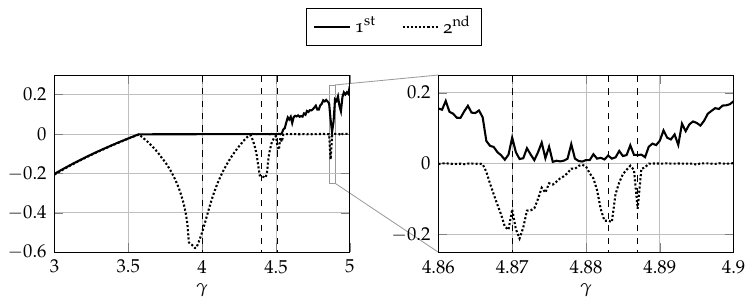}
\par
\includegraphics[trim=19.215pt 0 12.782pt 0]{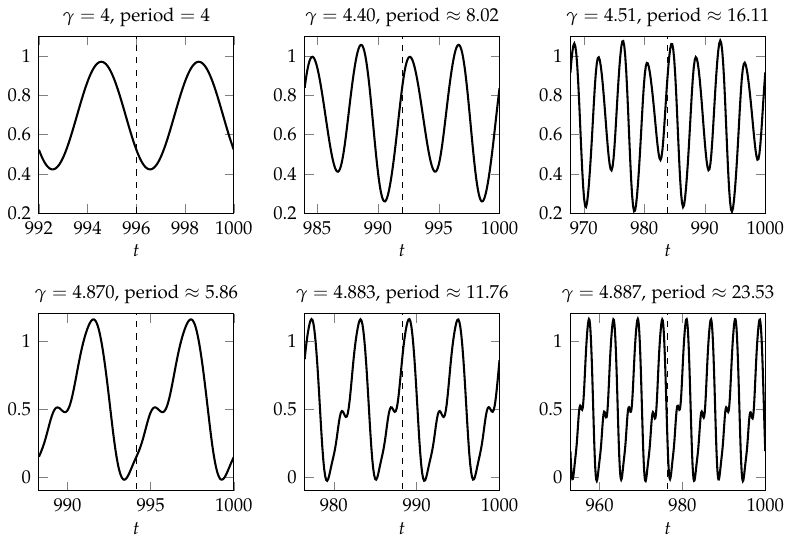}
\caption{Diagram of the first two dominant (in descending order) LEs (top row) and solutions (other rows) of the RE with quadratic nonlinearity \cref{quadRE} when varying $\gamma$, computed with $M_X=15$ and $T=1000$.
The solutions are computed via pseudospectral discretization with $M_X=15$, starting from a constant initial value of $0.2$.
The final time of $T=1000$ ensures sufficiently good convergence to the stable periodic solution.
For each solution, the last two periods are shown, separated by a vertical dashed line.
The values of $\gamma$ and of the period are given above the diagrams; the values of $\gamma$ are also marked by vertical dashed lines in the diagram of the LEs.}
\label{fig:quad-le-diagram-persol}
\end{center}
\end{figure}

\Cref{fig:quad-le-diagram-persol} (compare with \cite[Figure 2.3]{BredaDiekmannLiessiScarabel2016}) presents, on the top row, the diagram of the first two dominant LEs of \cref{quadRE} when varying $\gamma$, computed with $M_X=15$ and $T=1000$, following previous experimental considerations.
We can observe that the dominant LE is $0$ at the expected Hopf bifurcation ($\gamma=2+\frac{\pi}{2}$), after which one LE is always $0$ since periodic solutions appear.
The dominant nontrivial exponent reaches $0$ again for the expected period-doubling bifurcations at $\gamma\approx 4.32, 4.49, 4.53$, and it becomes positive for $\gamma\geq4.55$, indicating the insurgence of a chaotic regime, which is compatible with what was obtained in \cite{BredaDiekmannLiessiScarabel2016}.
Finally, for $\gamma \in [4.86,4.9]$ other \emph{stability islands} appear, corresponding to a branch of stable periodic solutions (appearing at $\gamma\approx4.8665$) and the corresponding cascade of period-doubling bifurcations (at $\gamma\approx4.8795,4.8860$) leading back to chaos (starting at $\gamma\approx4.8885$).
As an example, \cref{fig:quad-le-diagram-persol} (second and third row) shows some stable periodic solutions in the branches arising from the first and the second set of bifurcations.
Observe that indeed the period approximately doubles at each bifurcation.

\subsection{Egg cannibalism model}
\label{sec:cannibRE}

The second equation we consider is the egg cannibalism (toy) model from \cite{BredaDiekmannMasetVermiglio2013}:
\begin{equation}\label{cannibRE}
x(t) = \frac{\gamma}{2}\int_{a_{\text{mat}}}^{a_{\text{max}}} x(t-a)\E^{-x(t-a)} \D a,
\end{equation}
with $a_{\text{mat}}$ and $a_{\text{max}}$ being, respectively, the constant maturation and maximum ages.
Observe that \cref{cannibRE} corresponds to \cref{REprototype} with $\tau_1\coloneqq a_{\text{mat}}$, $\tau_2\coloneqq a_{\text{max}}$ and $f(x) \coloneqq \frac{\gamma}{2}x\E^{-x}$.
Also in this case, the equilibria and their stability properties are known, including the occurrence of a Hopf bifurcation for the nontrivial equilibrium at $\log(\gamma) = 1+\frac{\pi}{2}$, although here the periodic solutions are not explicitly known; again, the presence of period-doubling bifurcations is experimentally known \cite{BredaDiekmannGyllenbergScarabelVermiglio2016,BredaDiekmannLiessiScarabel2016,ScarabelDiekmannVermiglio2021,BredaDiekmannMasetVermiglio2013}.

\begin{figure}
\begin{center}
\includegraphics[trim=24.937pt 0 44.242pt 0]{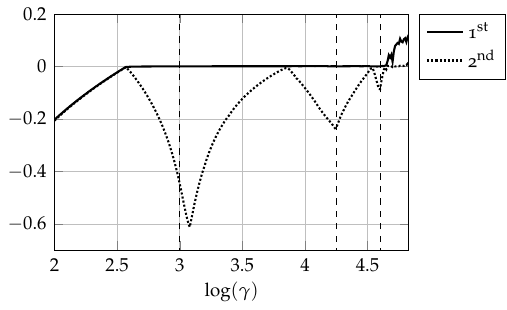}
\par
\includegraphics[trim=16.965pt 0 12.782pt 0]{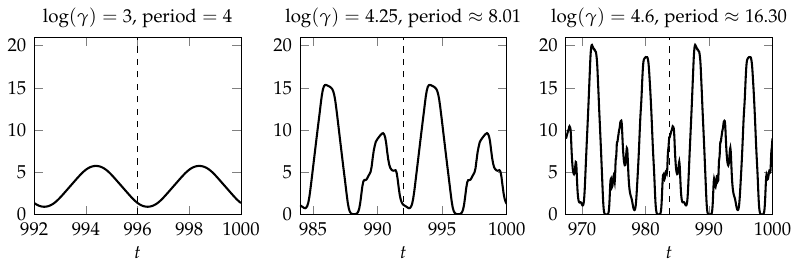}
\caption{Diagram of the first two dominant (in descending order) LEs (top row) and solutions (bottom row) of the egg cannibalism RE \cref{cannibRE} with $a_{\text{mat}}=1$ and $a_{\text{max}}=3$, when varying $\gamma$, computed with $M_X=15$ and $T=1000$.
The solutions are computed via pseudospectral discretization with $M_X=15$, starting from a constant initial value of $0.2$.
The final time of $T=1000$ ensures sufficiently good convergence to the stable periodic solution.
For each solution, the last two periods are shown, separated by a vertical dashed line.
The values of $\gamma$ and of the period are given above the diagrams; the values of $\gamma$ are also marked by vertical dashed lines in the diagram of LEs.}
\label{fig:cannib-le-diagram-persol}
\end{center}
\end{figure}

\Cref{fig:cannib-le-diagram-persol} (top row) presents the diagram of the first two dominant LEs of \cref{cannibRE} when varying $\gamma$, with $a_{\text{mat}}=1$ and $a_{\text{max}}=3$.
The numerical parameters are $M_X=15$, $T=1000$, a tolerance of $10^{-6}$ for \texttt{dqr}, and $\texttt{RelTol} = 10^{-6}$ and $\texttt{AbsTol} = 10^{-7}$ for \texttt{ode45}.
Similar to the previous example, the dominant exponent is $0$ at the Hopf bifurcation, and one exponent remains $0$ thereafter.
The dominant nontrivial exponent reaches $0$ again for the expected period-doubling bifurcations at $\log(\gamma)\approx 3.855$ (\cite{ScarabelDiekmannVermiglio2021} finds $\log(\gamma)\approx3.8777$ with $M_X=20$ and $\log(\gamma)\approx3.8763$ with $M_X=40$, setting MatCont's tolerances to $10^{-10}$ for Newton's method and $10^{-6}$ for the calculation of the test functions for bifurcation points) and $\log(\gamma)\approx 4.54$, and it becomes positive for $\log(\gamma)\geq4.66$, indicating chaos.
\Cref{fig:cannib-le-diagram-persol} (bottom row) shows some stable periodic solutions, confirming the (approximate) doubling of the period.

\subsection{Simplified logistic Daphnia}
\label{sec:logDaphnia}

The third and final equation is a simplified version of the Daphnia model with a logistic term as the growth of the resource, taken from \cite{BredaDiekmannMasetVermiglio2013}:
\begin{equation}\label{logDaphnia}
\left\{
\begin{aligned}
&b(t) = \beta S(t) \int_{a_{\text{mat}}}^{a_{\text{max}}} b(t-a) \D a, \\
&S'(t) = rS(t)\biggl(1-\frac{S(t)}{K}\biggr) -\gamma S(t) \int_{a_{\text{mat}}}^{a_{\text{max}}} b(t-a) \D a,
\end{aligned}
\right.
\end{equation}
where $b$ is the birth rate of the consumer population, $S$ is the density of the resource, $r$ and $K$ are, respectively, the growth rate and the carrying capacity of the resource, and $a_{\text{mat}}$ and $a_{\text{max}}$ have the same meaning as in \cref{sec:cannibRE}%
\footnote{%
Note that the second term on the right-hand side of the equation for $S$ in \cref{logDaphnia} can be rewritten as $\gamma b(t)/\beta$, rendering that equation an ODE.
This does not bring any simplification since the integral term is computed for the first equation anyway.
We prefer to keep the form of \cref{logDaphnia} because it comes from a more general class where the fertility and consumption rates are not constant, but, rather, are functions of the age and the size of the individuals (see (1) in \cite{BredaDiekmannMasetVermiglio2013}).
}%
.
\Cref{logDaphnia} corresponds to \cref{coupledprototype} with $\tau_1\coloneqq a_{\text{mat}}$, $\tau_2\coloneqq a_{\text{max}}$, $f_1(x) \coloneqq \beta x$, $f_2(x) \coloneqq -\gamma x$ and $g(y) \coloneqq ry(1-\frac{y}{K})$.
The equilibria are known, along with the stability properties of the trivial (zero and consumer-free) ones; in particular, the consumer-free equlibrium exchanges stability with the nontrivial equilibrium in a transcritical bifurcation at $\beta=(K(a_{\text{max}}-a_{\text{mat}}))^{-1}$; moreover, when varying $\beta$, the nontrivial equilibrium is experimentally known to undergo a Hopf bifurcation \cite{BredaDiekmannGyllenbergScarabelVermiglio2016,BredaDiekmannMasetVermiglio2013,BredaLiessi2020}.

\begin{figure}
\begin{center}
\includegraphics[trim=24.937pt 0 44.240pt 0]{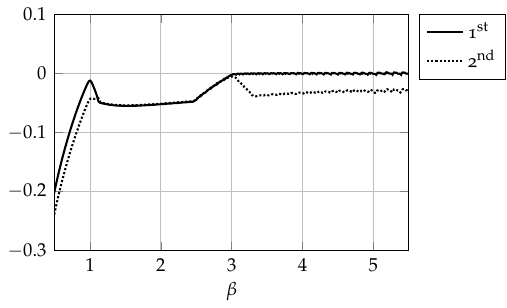}
\caption{Diagram of the first two dominant (in descending order) LEs of the simplified logistic Daphnia equation \cref{logDaphnia} with $a_{\text{mat}}=3$, $a_{\text{max}}=4$ and $r=K=\gamma=1$, when varying $\beta$, computed with $M_X=15$ and $T=1000$.}
\label{fig:logdaphnia-le-diagram}
\end{center}
\end{figure}

The diagram of the first two dominant LEs of \cref{logDaphnia} when varying $\beta$ is shown in \cref{fig:logdaphnia-le-diagram}.
The values of the other model parameters are $a_{\text{mat}}=3$, $a_{\text{max}}=4$ and $r=K=\gamma=1$.
The numerical parameters are $M_X=M_Y=15$, $T=1000$, a tolerance of $10^{-6}$ for \texttt{dqr}, and $\texttt{RelTol} = 10^{-4}$ and $\texttt{AbsTol} = 10^{-5}$ for \texttt{ode45}.
We can observe a spike at $\beta=1$ (albeit not touching the value $0$), corresponding to the known transcritical bifurcation, while the Hopf bifurcation seems to happen for $\gamma\approx3$ (\cite{BredaDiekmannGyllenbergScarabelVermiglio2016} finds $\gamma\approx3.0161$ with $M_X=10$ and MatCont's tolerances set to $10^{-10}$).
We continued the experiment for values of $\beta$ higher than those shown in \cref{fig:logdaphnia-le-diagram}, but we did not find other bifurcations.

Compared with the previous examples, the diagram seems less accurate (observe the spike corresponding to the transcritical bifurcation and the trivial LE after the Hopf bifurcation).
The explanation for this phenomenon is still unknown.
In this example, we have increased the tolerances for \texttt{ode45} in order to improve the computation times; however, for $\beta=1$ the dominant LE differs absolutely from the one computed with $\texttt{RelTol} = 10^{-6}$ and $\texttt{AbsTol} = 10^{-7}$ (as in the previous examples) by less than $10^{-7}$ and is thus still substantially far from $0$ (as a further example, for $\beta=5$ the absolute difference is larger but still less than $10^{-5}$).
Moreover, as $\beta$ increases, the periodic solution presents flat regions followed by spikes, which may suggest that the equation is becoming stiff; however, we tried replacing \texttt{ode45} with MATLAB's \texttt{ode23s}, implementing a modified Rosenbrock formula of order $2$ for stiff ODEs \cite{ShampineReichelt1997}, with no substantial changes.

\section{Concluding remarks}
\label{sec:concluding}

In this work, we have provided the first method, to our knowledge, for computing LEs of REs and coupled equations.
The proposed method appears to be effective when applied to examples with known properties; however, since the nature of our study has been purely computational, further investigation into the method's convergence properties is required.
As far as efficiency is concerned, LEs are notoriously expensive to compute \cite{BredaVanVleck2014}, and that is true also in this case; the computational cost depends linearly on $T$, while its relation with $M_{X}$ and $M_{Y}$ is currently unclear from the experiments, even though it is expected to be of order $4$, according to \cref{IVPGAMMA}.

The next step in our research is to tackle the problem following the approach of \cite{BredaVanVleck2014}; as recalled in \cref{sec:intro}, this involves the reformulation of the equation (RE or coupled equation, in our case) in a Hilbert space and the rigorous definition of LEs and the DQR in the new setting.
The technique of \cite{BredaDellaSchiava2018} and of the present work, however, is a rather general approach which can also be applied to other kinds of equations (e.g., certain classes of partial differential equations of interest for mathematical biology \cite{AbiaAnguloLopezMarcos2022}), as anticipated in \cite{BredaDellaSchiava2018} and according to the \emph{pragmatic} point of view discussed in \cite{BredaDiekmannLiessiScarabel2016}.

\section*{Acknowledgments}

The authors are members of INdAM research group GNCS and of UMI research group ``Mo\-del\-li\-sti\-ca socio-epidemiologica''.
This work was partially supported by the Italian Ministry of University and Research (MUR) through the PRIN 2020 project (No.\ 2020JLWP23) ``Integrated Mathematical Approaches to Socio-Epidemiological Dynamics'' (CUP: E15F21005420006), and by the GNCS 2023 project ``Sistemi dinamici e modelli di evoluzione: tecniche funzionali, analisi qualitativa e metodi numerici'' (CUP: E53C22001930001).

{\sloppy
\printbibliography
\par}

\end{document}